\documentclass[12pt,reqno]{amsart}
\usepackage{indentfirst, amssymb,amsmath,amsthm}
\newtheorem{theo}{Theorem}[section]

\newtheorem{lem}{Lemma}[section]

\newtheorem{cor}{Corollary}[section]

\newcommand{\be}{\begin{equation}}
\newcommand{\ee}{\end{equation}}
\newcommand{\beas}{\begin{eqnarray*}}
\newcommand{\eeas}{\end{eqnarray*}}
\newcommand{\bea}{\begin{eqnarray}}
\newcommand{\eea}{\end{eqnarray}}

\numberwithin{equation}{section}

\begin{document}

\setlength{\unitlength}{1mm} \baselineskip .45cm
%\large
\setcounter{page}{1}
\pagenumbering{arabic}
\title[ $h$-almost Ricci-Yamabe solitons in paracontact geometry]
{$h$-almost Ricci-Yamabe solitons in paracontact geometry}

\author[A. Sardar, U. C. De and C. \"Ozg\"ur]{Arpan Sardar*, Uday Chand De and Cihan \"Ozg\"ur}

\address{Arpan Sardar\newline  Department of Mathematics,\newline University of Kalyani,\newline Kalyani 741235, West Bengal, India}
	\email {arpansardar51@gmail.com}

\address{Uday Chand De \newline  Department of Pure Mathematics,\newline University of Calcutta,
\newline 35, Ballygaunge Circular Road,\newline Kolkata 700019, West Bengal, India}
 \email {uc$\_$de@yahoo.com}

\address{Cihan \"Ozg\"ur \newline  Department of Mathematics,\newline Izmir Democracy University,
\newline Karabaglar, 35140,Izmir/Turkey}
 \email {cihan.ozgur@idu.edu.tr}

\begin{abstract}
In this article, we classify $h$-almost Ricci-Yamabe solitons in paracontact geometry. In particular, we characterize para-Kenmotsu manifolds satisfying $h$-almost Ricci-Yamabe solitons and 3-dimensional para-Kenmotsu manifolds obeying  $h$-almost gradient Ricci-Yamabe solitons. Next, we classify para-Sasakian manifolds admitting $h$-almost Ricci-Yamabe solitons and  $h$-almost gradient Ricci-Yamabe solitons. Besides these, we investigate $h$-almost Ricci-Yamabe solitons and $h$-almost gradient Ricci-Yamabe solitons in para-cosymplectic manifolds. Finally, we construct two examples to illustrate our results.
\end{abstract}

\maketitle
\footnotetext{\\
\\
*Arpan Sardar is financially supported by UGC, Ref. ID. 4603/(CSIR-UGCNETJUNE2019).\\
}

\section{\textsf{Introduction}} 

In a semi-Riemannian manifold the Ricci-Yamabe soliton  is defined by
\begin{equation}\label{1.1}
\pounds_V g + (2\lambda-\beta r)g + 2\alpha S = 0,
\end{equation}
$\pounds$ denotes the Lie-derivative, $S$ is the Ricci tensor, $r$ being the scalar curvature and $\lambda, \alpha,\beta \in \mathbb{R}$.  Ricci-Yamabe solitons are the special solutions of the Ricci-Yamabe flow
\begin{equation}
\frac{\partial g}{\partial t} = -2\alpha S + \beta r g,
\end{equation}
which was introduced by Guler and Crasmareanu\cite{guler}. Equation (\ref{1.1}) is called almost Ricci-Yamabe soliton provided $\lambda$ is a smooth function.\\

In particular,  for $\alpha = 1$ and $\beta =0$, (\ref{1.1}) implies
\begin{equation}
\pounds_V g + 2\alpha S + 2\lambda g = 0,
\end{equation}
which is a Ricci soliton for $\lambda \in \mathbb{R}$. Thus almost Ricci-Yamabe solitons ( Ricci-Yamabe solitons) are the natural generalization of almost Ricci solitons (Ricci solitons). Several generalization of  Ricci solitons are almost Ricci solitons(\cite{des1}, \cite{des2}, \cite{wan}, \cite{wan2}, \cite{wan1}, \cite{w7}), $\eta$-Ricci solitons (\cite{bla1}, \cite{bla2}, \cite{bla3}, \cite{de}, \cite{sar}, \cite{sar1}), $\ast$-Ricci solitons(\cite{dai}, \cite{ham}, \cite{kai}, \cite{ven}, \cite{wan3}) and many others.\\

Recently, Gomes et al.\cite{gom}  extended the concept of almost Ricci soliton to h-almost Ricci soliton on a complete Riemannian manifold  by
\begin{equation}\label{1.2}
\frac{h}{2}\pounds_V g + \lambda g + S = 0,
\end{equation}
where $h: M \rightarrow \mathbb{R}$ is a smooth function. Specifically, a Ricci soliton is the 1-almost Ricci soliton endowed with constant $\lambda$.\\

Now we introduce the new type of solitons named $h$-almost Ricci-Yamabe solitons (briefly, h-ARYS) which are the extended notion of almost Ricci-Yamabe solitons, which are given by
\begin{equation}\label{1.3}
\frac{h}{2}\pounds_V g + \alpha S + (\lambda-\frac{\beta}{2}r) g = 0,
\end{equation}
where $h$ is a smooth function.\\

If $V$ is a gradient of $f$ on the manifold, then the foregoing concept is called $h$-almost gradient Ricci-Yamabe soliton (briefly, h-AGRYS) and (\ref{1.3}) takes the form
\begin{equation}\label{1.4}
h\nabla^2 f + (\lambda -\frac{\beta}{2}r)g + \alpha S = 0.
\end{equation}
An $h$-AGRYS is named $h$-gradient Ricci-Yamabe soliton if $\lambda$ is constant.\\

An h-ARYS (or h-AGRYS) turns into\\

(i) $h$-almost Ricci soliton (or $h$-almost gradient Ricci soliton) if $\beta = 0$ and $\alpha = 1$,\\

(ii) $h$-almost Yamabe soliton (or $h$-almost gradient Yamabe soliton) if $\beta = 1$ and $\alpha = 0$,\\

(iii) $h$-almost Einstein soliton (or $h$-almost gradient Einstein soliton) if $\beta = -1$ and $\alpha = 1$.\\

The $h$-ARYS ( or  $h$-AGRYS ) is called proper if $\alpha \neq 0,1$.\\

Recently, in (\cite{sar2}, \cite{sar3}), the first author and Sarkar studied Ricci-Yamabe solitons in Kenmotsu 3-manifolds and generalized Sasakian space forms. Also, Sing and Khatri\cite{kha} studied Ricci-Yamabe solitons in perfect fluid spacetime.\\\\

The current article is structured as:

After the introduction, required preliminaries have been mentioned in Section 2. In Section 3, we investigate $h$-ARYS and $h$-AGRYS in para-Kenmotsu manifolds. Next, we classify para-Sasakian manifolds admitting $h$-ARYS and  $h$-AGRYS in Section 5. Besides these, we investigate $h$-ARYS and $h$-AGRYS in para-cosymplectic manifolds in Section 7. Finally, we construct two  examples to illustrate our results.

\section{\textsf{Preliminaries}}

 An almost paracontact structure on a manifold $M^{2n+1}$ consists of a (1,1)-tensor field $\phi$, a vector field $\zeta$ and a one-form $\eta$ obeying the subsequent conditions:
\begin{equation}\label{2.1}
\phi^2 = I - \eta\otimes \zeta, \hspace{.4cm} \eta(\zeta) = 1
\end{equation}
and the tensor field $\phi$ induces an almost paracomplex structure on each fibre of $\mathcal{D} = ker(\eta)$, that is, the $\pm 1$-eigendistributions, $\mathcal{D}^\pm = \mathcal{D}_\phi (\pm 1)$ of $\phi$ have equal dimension $n$.

The manifold $M$ with an almost paracontact structure is named an almost paracontact manifold. From the definition it can be established that $\phi \zeta = 0$, $\eta\circ \phi = 0$ and rank of $\phi$ is $2n$. If the Nijenhuis tensor vanishes identically, then this manifold is said to be normal. $M$ is named an almost paracontact metric manifold if there exists a semi-Riemannian metric $g$ such that
\begin{equation}\label{2.2}
g(\phi Z_1, \phi Z_2) = -g(Z_1,Z_2) + \eta(Z_1)\eta(Z_2)
\end{equation}
for all $Z_1,Z_2 \in \chi(M)$.\\

$(M,\phi,\zeta,\eta,g)$ is named a paracontact metric manifold if $d\eta(Z_1,Z_2) = g(Z_1,\phi Z_2) = \Phi(Z_1,Z_2)$, $\Phi$ being the fundamental 2-form of $M$.\\

 An almost paracontact metric manifold $M^{2n+1}$, with a structure $(\phi,\zeta,\eta,g)$ is said to be an almost $\gamma$-paracosymplectic manifold, if
\begin{equation}\label{2.3}
d\eta = 0, d\Phi = 2\gamma \eta\wedge \Phi,
\end{equation}
where $\gamma$ is a constant or function on $M$. If we put $\gamma = 1$ in (\ref{2.3}), we acquire almost para-Kenmotsu manifolds. In \cite{erk}, the para-Kenmotsu manifold satisfies
\begin{equation}\label{2.4}
R(Z_1,Z_2)\zeta = \eta(Z_1)\;Z_2 - \eta(Z_2)\;Z_1,
\end{equation}
\begin{equation}\label{2.5}
R(Z_1,\zeta)\;Z_2= g(Z_1,\;Z_2)\zeta -\eta(Z_2)\;Z_1,
\end{equation}
\begin{equation}\label{2.6}
R(\zeta,Z_1)Z_2 = -g(Z_1,Z_2)\zeta +\eta(Z_2)Z_1,
\end{equation}
\begin{equation}\label{2.7}
\eta(R(Z_1,Z_2)Z_3) = - g(Z_2,Z_3)\eta(Z_1)+g(Z_1,Z_3)\eta(Z_2) ,
\end{equation}
\begin{equation}\label{2.8}
(\nabla_{Z_1} \phi)Z_2 = g(\phi Z_1,Z_2)\zeta -\eta(Z_2)\phi Z_1,
\end{equation}
\begin{equation}\label{2.9}
\nabla_{Z_1} \zeta = Z_1 - \eta(Z_1)\zeta,
\end{equation}
\begin{equation}\label{2.10}
S(Z_1,\zeta) = -2n \eta(Z_1).
\end{equation}

\begin{lem}(\cite{erk})
In a para-Kenmotsu manifold $M^3$, we have
\begin{equation}\label{2.11}
\zeta r = -2(r+6).
\end{equation}
\end{lem}
Also in a $M^3$, we have
\begin{equation}\label{2.12}
QZ_1 = (\frac{r}{2}+1)Z_1 -(\frac{r}{2}+3)\eta(Z_1)\zeta,
\end{equation}
which implies
\begin{equation}\label{2.13}
S(Z_1, Z_2) =-(\frac{r}{2}+3)\eta(Z_1)\eta(Z_2)+ (\frac{r}{2}+1)g(Z_1, Z_2) ,
\end{equation}
 $Q$ indicates the Ricci operator defined by $S(Z_1,Z_2) = g(QZ_1,Z_2)$.

\vspace{.6cm}

{\section{\textsf{$h$-ARYS on para-Kenmotsu manifolds}}}

We assume that the manifold $M^{2n+1}$ admits an $h$-ARYS $(g,\zeta, \lambda, \alpha, \beta)$. Then from (\ref{1.3}), we get
\begin{equation}\label{3.1}
\frac{h}{2}(\pounds_{\zeta} g)(Z_1,Z_2) + \alpha S(Z_1,Z_2) + (\lambda -\frac{\beta}{2}r)g(Z_1,Z_2) = 0,
\end{equation}
which implies
\begin{equation}\label{3.2}
\frac{h}{2}[g(\nabla_{Z_1} \zeta ,Z_2) + g(Z_1,\nabla_{Z_2} \zeta)] + \alpha S(Z_1,Z_2) + (\lambda -\frac{\beta}{2}r)g(Z_1,Z_2) = 0,
\end{equation}
Using (\ref{2.9}) in (\ref{3.2}), we infer
\begin{equation}\label{3.3}
\alpha S(Z_1,Z_2) =  h\eta(Z_1)\eta(Z_2)-(h+\lambda -\frac{\beta}{2}r)g(Z_1,Z_2) .
\end{equation}
Putting $Z_1 = Z_2 = \zeta$ in the foregoing equation entails that
\begin{equation}\label{3.4}
\frac{\beta}{2}r = \lambda - 2n\alpha.
\end{equation}
Equations (\ref{3.3}) and (\ref{3.4}) together give
\begin{equation}
\alpha S(Z_1,Z_2) = -(h+2n\alpha)g(Z_1,Z_2) + h\eta(Z_1)\eta(Z_2),
\end{equation}
which is an $\eta$-Einstein manifold. Hence we have:\\

\begin{theo}
If a $M^{2n+1}$ admits a proper $h$-ARYS, then the manifold becomes an $\eta$-Einstein manifold.
\end{theo}

\vspace{1cm}

Let $M^3$ admit an $h$-AGRYS. Then (\ref{1.4}) implies
\begin{equation}\label{4.1}
h\nabla_{Z_1} Df = -\alpha Q Z_1 - (\lambda -\frac{\beta}{2}r)Z_1.
\end{equation}
Taking covariant derivative of (\ref{4.1}), we get
\begin{eqnarray}\label{4.2}
h\nabla_{Z_2} \nabla_{Z_1} Df &=& \frac{1}{h}(Z_2 h)[\alpha Q Z_1 + (\lambda-\frac{\beta}{2}r)Z_1] - \alpha \nabla_{Z_2} QZ_1\\ \nonumber
&& -(Z_2\lambda)Z_1 -(\lambda-\frac{\beta}{2}r)\nabla_{Z_2} Z_1 + \frac{\beta}{2}(Z_2 r)Z_1.
\end{eqnarray}
Interchanging $Z_1$ and $Z_2$ in (\ref{4.2}) entails that
\begin{eqnarray}\label{4.3}
h\nabla_{Z_1} \nabla_{Z_2} Df &=& \frac{1}{h}(Z_1 h)[\alpha QZ_2 + (\lambda-\frac{\beta}{2}r)Z_2] - \alpha \nabla_{Z_1} QZ_2\\ \nonumber
&& -(Z_1\lambda)Z_2 -(\lambda-\frac{\beta}{2}r)\nabla_{Z_1} Z_2 + \frac{\beta}{2}(Z_1 r)Z_2.
\end{eqnarray}
Equation (\ref{4.1}) implies
\begin{equation}\label{4.4}
h\nabla_{[Z_1, Z_2]} Df = -\alpha Q([Z_1, Z_2]) - (\lambda -\frac{\beta}{2}r)([Z_1, Z_2]).
\end{equation}
Equations (\ref{4.2}), (\ref{4.3}) and (\ref{4.4}) reveal that
\begin{eqnarray}\label{4.5}
hR(Z_1,Z_2)Df &=& \frac{1}{h}(Z_1 h)[\alpha QZ_2 + (\lambda-\frac{\beta}{2}r)Z_2]\\ \nonumber
&& - \frac{1}{h}(Z_2 h)[\alpha QZ_1 + (\lambda-\frac{\beta}{2}r)Z_1]\\ \nonumber
&& -\alpha[(\nabla_{Z_1} Q)Z_2 - (\nabla_{Z_2} Q)Z_1]\\ \nonumber
&& + \frac{\beta}{2}[(Z_1 r)Z_2-(Z_2 r)Z_1] - [(Z_1 \lambda)Z_2 - (Z_2 \lambda)Z_1].
\end{eqnarray}
Equation (\ref{2.12}) implies
\begin{eqnarray}\label{4.6}
(\nabla_{Z_1} Q)Z_2 &=& \frac{Z_1 r}{2}[Z_2-\eta(Z_2)\zeta]\\ \nonumber
&& -(3+\frac{r}{2})[g(Z_1,Z_2)\zeta -2\eta(Z_1)\eta(Z_2)\zeta + \eta(Z_2)Z_1].
\end{eqnarray}
Using (\ref{4.6}) in (\ref{4.5}), we get
\begin{eqnarray}\label{4.7}
hR(Z_1,Z_2)Df &=&  \frac{1}{h}(Z_1 h)[\alpha QZ_2 + (\lambda-\frac{\beta}{2}r)Z_2]\\ \nonumber
&& - \frac{1}{h}(Z_2 h)[\alpha QZ_1 + (\lambda-\frac{\beta}{2}r)Z_1]\\ \nonumber
&& -\alpha \frac{(Z_1 r)}{2}[Z_2 -\eta(Z_2)\zeta] + \alpha \frac{(Z_2 r)}{2}[Z_1 -\eta(Z_1)\zeta]\\ \nonumber
&& + \alpha(3+\frac{r}{2})[\eta(Z_2)Z_1 -\eta(Z_1)Z_2] - [(Z_1\lambda)Z_2-(Z_2\lambda )Z_1]\\\nonumber
&& + \frac{\beta}{2}[(Z_1 r)Z_2 - (Z_2 r)Z_1].
\end{eqnarray}
If we take $h$ = constant, then (\ref{4.7}) implies
\begin{eqnarray}\label{4.8}
hR(Z_1,Z_2)Df &=& -\alpha \frac{(Z_1 r)}{2}[Z_2 -\eta(Z_2)\zeta] + \alpha \frac{(Z_2 r)}{2}[Z_1 -\eta(Z_1)\zeta]\\ \nonumber
&& + \alpha(3+\frac{r}{2})[\eta(Z_2)Z_1 -\eta(Z_1)Z_2] - [(Z_1\lambda)Z_2(Z_2\lambda )Z_1]\\\nonumber
&& + \frac{\beta}{2}[(Z_1 r)Z_2 - (Z_2 r)Z_1].
\end{eqnarray}
Contracting (\ref{4.8}), we infer
\begin{equation}\label{4.9}
hS(Z_2,Df) = (\frac{\alpha}{2}-\beta)Z_2 r + 2(Z_2 \lambda).
\end{equation}
Replacing $Z_1$ by $Df$ in (\ref{2.13}) and comparing with (\ref{4.9}), we get
\begin{equation}\label{4.10}
h[(1+\frac{r}{2})Z_2 f -(3+\frac{r}{2})(\zeta f)\eta(Z_2)] = (\frac{\alpha}{2}-\beta)Z_2 r + 2(Z_2\lambda).
\end{equation}
Putting $Z_2=\zeta$ in (\ref{4.10}) entails that
\begin{equation}\label{4.11}
h(\zeta f) = (\frac{\alpha}{2}-\beta)(r+6) - (\zeta \lambda).
\end{equation}
Taking inner product of (\ref{4.8}) with $\zeta$ gives
\begin{eqnarray}\label{4.12}
h[\eta(Z_2)Z_1f - \eta(Z_1)Z_2f] &=& -[(Z_1\lambda)\eta(Z_2)-(Z_2\lambda)\eta(Z_1)]\\ \nonumber
&& + \frac{\beta}{2}[(Z_1r)\eta(Z_2)-(Z_2r)\eta(Z_1)].
\end{eqnarray}
Setting $Z_2 =\zeta$ in (\ref{4.12}) and using (\ref{4.11}), we get
\begin{equation}\label{4.13}
h(Z_1f) = \frac{\beta}{2}(Z_1 r) + \frac{\alpha}{2}(r+6)\eta(Z_1) - (Z_1\lambda).
\end{equation}
Let us assume that the scalar curvature $r$ = constant. Then from (\ref{2.11}) we get $r = -6$. Therefore the above equation implies
\begin{equation}\label{4.14}
h(Z_1 f) = -(Z_1 \lambda),
\end{equation}
which implies
\begin{equation}\label{4.15}
h (Df) =  - (D\lambda).
\end{equation}
Using (\ref{4.15}) in (\ref{4.1}) reveals that
\begin{equation}\nonumber
-\nabla_{Z_1} D\lambda = -\alpha QZ_1 - (\lambda -\frac{\beta}{2}r)Z_1,
\end{equation}
which shows that it is an almost gradient Ricci-Yamabe soliton whose soliton function is $-\lambda$. Hence we have:\\

\begin{theo}
If a $M^3$ with a constant scalar curvature admits an $h$-ARYS, then the soliton becomes an almost gradient Ricci-Yamabe soliton whose soliton function is -$\lambda$.
\end{theo}

\vspace{.9cm}

{\section{\textsf{para-Sasakian manifolds}}}

A para-Sasakian manifold is a  normal paracontact metric manifold. It is to be noted that the para-Sasakian manifold implies the $K$-paracontact manifold and conversely (only in 3 dimensions). In a para-Sasakian manifold  the following relations hold:
\begin{equation}\label{6.1}
R(Z_1,Z_2)\zeta = \eta(Z_1)Z_2 - \eta(Z_2)Z_1,
\end{equation}
\begin{equation}\label{6.2}
(\nabla_{Z_1} \phi)Z_2 = -g(Z_1,Z_2)\zeta + \eta(Z_2)Z_1,
\end{equation}
\begin{equation}\label{6.3}
\nabla_{Z_1} \zeta = -\phi Z_1,
\end{equation}
\begin{equation}\label{6.4}
R(\zeta,Z_1)Z_2 = -g(Z_1,Z_2)\zeta + \eta(Z_2)Z_1,
\end{equation}
\begin{equation}\label{6.5}
S(Z_1,\zeta) = -2n\eta(Z_1).
\end{equation}
In a 3-dimensional semi-Riemannian manifold the curvature tensor $R$ is of the form
\begin{eqnarray}\label{6.6}
R(Z_1,Z_2)Z_3 &=& g(Z_2,Z_3)QZ_1 - g(Z_1,Z_3)QZ_2 + S(Z_2,Z_3)Z_1 \\ \nonumber
&&- S(Z_1,Z_3)Z_2 - \frac{r}{2}[g(Z_2,Z_3)Z_1 - g(Z_1,Z_3)Z_2].
\end{eqnarray}

Equation (\ref{6.6}) implies
\begin{equation}\label{6.7}
QZ_1 = (\frac{r}{2}+1)Z_1-(\frac{r}{2}+3)\eta(Z_1)\zeta,
\end{equation}
which implies
\begin{equation}\label{6.8}
S(Z_1,Z_2) = (\frac{r}{2}+1)g(Z_1,Z_2) -(\frac{r}{2}+3)\eta(Z_1)\eta(Z_2).
\end{equation}

\begin{lem}(\cite{erk})
For a para-Sasakian manifold $M^3$, we have
\begin{equation}\label{6.9}
\zeta r = 0.
\end{equation}
\end{lem}

\vspace{.6cm}

{\section{\textsf{$h$-ARYS on  para-Sasakian manifolds}}}

Let us assume that a  para-Sasakian manifold $M^{2n+1}$ admit an $h$-ARYS $(g, \zeta, \lambda, \alpha, \beta)$. Then equation (\ref{1.3}) implies
\begin{equation}\label{7.1}
\frac{h}{2}(\pounds_{\zeta} g)(Z_1,Z_2) + \alpha S(Z_1,Z_2) + (\lambda-\frac{\beta}{2}r)g(Z_1,Z_2) = 0,
\end{equation}
which implies
\begin{equation}\label{7.2}
\frac{h}{2}[g(\nabla_{Z_1} {\zeta} ,Z_2) + g(Z_1,\nabla_{Z_2} {\zeta})] + \alpha S(Z_1,Z_2) + (\lambda-\frac{\beta}{2}r)g(Z_1,Z_2) = 0.
\end{equation}
Using (\ref{6.3}) in (\ref{7.2}) entails that
\begin{equation}\label{7.3}
\alpha S(Z_1,Z_2) = (\frac{\beta}{2}r - \lambda)g(Z_1,Z_2).
\end{equation}
Putting $Z_1 = Z_2 = \zeta$ in (\ref{7.3}), we get
\begin{equation}\label{7.4}
\beta r = 2\lambda - 4n\alpha.
\end{equation}
 Hence from (\ref{7.3}), we infer
\begin{equation}\nonumber
S(Z_1,Z_2) = -2n g(Z_1,Z_2),
\end{equation}
since for proper $h$-ARYS ($\alpha \neq 0$). Hence it is an Einstein manifold. Therefore we state:\\

\begin{theo}
If $M^{2n+1}$ admits a proper $h$-ARYS, then the manifold becomes an Einstein manifold.
\end{theo}

If we take $\alpha = 1$ and $\beta = 0$, then (\ref{7.4}) implies $\lambda = 2n$. Hence we get:
\begin{cor}
If a $M^{2n+1}$ admits a proper $h$-almost Ricci soliton, then the soliton is expanding.
\end{cor}

\vspace{15cm}

Suppose that a $M^3$ admit an h-AGRYS. Then equation (\ref{1.4}) implies
\begin{equation}\label{8.1}
h\nabla_{Z_1} Df = -\alpha QZ_1 - (\lambda -\frac{\beta}{2}r)Z_1.
\end{equation}
Using (\ref{6.7}) in the above equation entails that
\begin{eqnarray}\label{8.2}
h\nabla_{Z_1} Df &=& -[\frac{(\alpha-\beta)}{2}r + \alpha + \lambda]Z_1\\ \nonumber
&& + \alpha(\frac{r}{2}+3)\eta(Z_1)\zeta.
\end{eqnarray}
Taking covariant differentiation of (\ref{8.2}), we get
\begin{eqnarray}\label{8.3}
h\nabla_{Z_2}\nabla_{Z_1} Df &=& \frac{1}{h}(Z_2 h)[\lbrace  \frac{(\alpha-\beta)}{2}r + \alpha + \lambda \rbrace Z_1 - \alpha(\frac{r}{2}+3)\eta(Z_1)\zeta ]\\ \nonumber
&&-  [\frac{(\alpha-\beta)}{2}Z_2 r + Z_2 \lambda]Z_1 -[\frac{(\alpha-\beta)}{2}r + \alpha + \lambda]\nabla_{Z_2} Z_1 \\ \nonumber
&&+ \frac{\alpha}{2}(Z_2 r)\eta(Z_1)\zeta + \alpha(\frac{r}{2}+3)[(\nabla_{Z_2} \eta(Z_1))\zeta - \eta(Z_1)\phi Z_2].
\end{eqnarray}
Swapping $Z_1$ and $Z_2$ in (\ref{8.3}), we infer
\begin{eqnarray}\label{8.4}
h\nabla_{Z_1}\nabla_{Z_2} Df &=& \frac{1}{h}(Z_1 h)[\lbrace  \frac{(\alpha-\beta)}{2}r + \alpha + \lambda \rbrace Z_2 - \alpha(\frac{r}{2}+3)\eta(Z_2)\zeta ]\\ \nonumber
&&-  [\frac{(\alpha-\beta)}{2}Z_1 r + Z_1 \lambda]Z_2 -[\frac{(\alpha-\beta)}{2}r + \alpha + \lambda]\nabla_{Z_1} Z_2 \\ \nonumber
&&+ \frac{\alpha}{2}(Z_1 r)\eta(Z_2)\zeta + \alpha(\frac{r}{2}+3)[(\nabla_{Z_1} \eta(Z_2))\zeta - \eta(Z_2)\phi Z_1].
\end{eqnarray}
Equation (\ref{8.2}) implies
\begin{eqnarray}\label{8.5}
h\nabla_{[Z_1,Z_2]} Df &=& -[\frac{(\alpha-\beta)}{2}r + \alpha + \lambda]([Z_1,Z_2])\\ \nonumber
&& + \alpha(\frac{r}{2}+3)\eta([Z_1,Z_2])\zeta.
\end{eqnarray}
With the help of (\ref{8.3})-(\ref{8.5}), we get
\begin{eqnarray}\label{8.6}
h R(Z_1,Z_2)Df &=& \frac{1}{h}(Z_1 h)[\lbrace  \frac{(\alpha-\beta)}{2}r + \alpha + \lambda \rbrace Z_2 - \alpha(\frac{r}{2}+3)\eta(Z_2)\zeta ]\\ \nonumber
&& - \frac{1}{h}(Z_2 h)[\lbrace  \frac{(\alpha-\beta)}{2}r + \alpha + \lambda \rbrace Z_1 - \alpha(\frac{r}{2}+3)\eta(Z_1)\zeta ]\\ \nonumber
&& -  [\frac{(\alpha-\beta)}{2}Z_1 r + Z_1\lambda]Z_2 +   [\frac{(\alpha-\beta)}{2}Z_2 r + Z_2\lambda]Z_1\\ \nonumber
&& + \frac{\alpha}{2}[(Z_1r)\eta(Z_2)\zeta - (Z_2 r)\eta(Z_1)\zeta]\\ \nonumber
&& + \alpha(\frac{r}{2}+3)[2g(Z_1,\phi Z_2)\zeta -\eta(Z_2)\phi Z_1 + \eta(Z_1)\phi Z_2].
\end{eqnarray}
Contracting the foregoing equation entails that
\begin{eqnarray}\label{8.9}
hS(Z_1,Df) &=& -\frac{1}{h}(Z_2h)[2\lbrace \frac{(\alpha-\beta)}{2}r + \alpha + \lambda \rbrace -\alpha(\frac{r}{2}+3)]\\ \nonumber
&& -\frac{\alpha}{h}(\frac{r}{2}+3)(\zeta h)\eta(Z_2) + (\frac{\alpha}{2}-\beta)(Z_2 r) + 2(Z_2 \lambda).
\end{eqnarray}
Replacing $Z_1$ by $Df$ in (\ref{6.8}) and likening with the above equation, we get
\begin{eqnarray}\label{8.10}
h[(\frac{r}{2}+1)Z_2 f -(\frac{r}{2}+3)(\zeta f) \eta(Z_2)] &=& -\frac{1}{h}(Z_2 h)[2\lbrace \frac{(\alpha-\beta)}{2}r + \alpha + \lambda \rbrace -\alpha(\frac{r}{2}+3)]\\ \nonumber
&& -\frac{\alpha}{h}(\frac{r}{2}+3)(\zeta h)\eta(Z_2) + (\frac{\alpha}{2}-\beta)(Z_2 r) + 2(Z_2 \lambda).
\end{eqnarray}
Setting $Z_2 = \zeta$ in (\ref{8.10}) reveals that
\begin{eqnarray}\label{8.11}
h(\zeta f) &=& \frac{1}{h}[\frac{(\alpha-\beta)}{2}r + \alpha + \lambda ](\zeta h) - (\zeta \lambda).
\end{eqnarray}
Taking inner product of (\ref{8.6}) with $\zeta$, we get
\begin{eqnarray}\label{8.12}
h[\eta(Z_2)Z_1 f - \eta(Z_1)Z_2 f] &=&\frac{1}{h}(Z_1 h)[\lbrace  \frac{(\alpha-\beta)}{2}r + \alpha + \lambda \rbrace  - \alpha(\frac{r}{2}+3)]\eta(Z_2)\\ \nonumber
&& - \frac{1}{h}(Z_2 h)[\lbrace  \frac{(\alpha-\beta)}{2}r + \alpha + \lambda \rbrace  - \alpha(\frac{r}{2}+3)]\eta(Z_1)\\ \nonumber
&& -[\frac{(\alpha-\beta)}{2}Z_1 r + Z_1\lambda]\eta(Z_2)\\ \nonumber
&& + [\frac{(\alpha-\beta)}{2}Z_2 r + Z_2\lambda]\eta(Z_1)\\ \nonumber
&&+ \frac{\alpha}{2}[(Z_1 r)\eta(Z_2)-(Z_2 r)\eta(Z_1)] + 2\alpha (\frac{r}{2}+3)g(Z_1,\phi Z_2).
\end{eqnarray}
Substituting $Z_1$ by $\phi Z_1$ and $Z_2$ by $\phi Z_2$ in (\ref{8.12}) gives
\begin{equation}
\alpha(r+6)g(\phi Z_1,Z_2) = 0.
\end{equation}
Since for proper  $\alpha \neq 0$, then the above equation implies $ r = -6$. Therefore from (\ref{6.8}), we get
\begin{equation}\label{8.14}
S(Z_1,Z_2) = -2 g(Z_1,Z_2),
\end{equation}
which is an Einstein manifold. In view of (\ref{6.6}) and (\ref{8.14}) reveals that
\begin{equation}
R(Z_1,Z_2)Z_3 = -[g(Z_2,Z_3)Z_1 -g(Z_1,Z_3)Z_2],
\end{equation}
which represents, it is a space of constant sectional curvature -1. Hence we have:\\

\begin{theo}
If a $M^3$ admits a proper $h$-ARYS, then the manifold is locally isometric to  $\mathbb{H}^3(1)$.
\end{theo}

\vspace{.9cm}

{\section{\textsf{para-cosymplectic manifolds}}}

An almost paracontact metric manifold $M^{2n+1}$ with a structure $(\phi,\zeta,\eta,g)$ is named an almost $\gamma$-paracosymplectic manifold if
\begin{equation}
d\eta = 0, d\Phi = 2\gamma \wedge \Phi.
\end{equation}

Specifically, if $\gamma = 0$, we acquire almost paracosymplectic manifolds. The manifold is called para-cosymplectic if it is normal. We refer (\cite{dak}, \cite{erk1}) for more details. Any paracosymplectic manifold satisfies
\begin{equation}\label{9.1}
R(Z_1,Z_2)\zeta = 0,
\end{equation}
\begin{equation}\label{9.2}
(\nabla_{Z_1} \phi )Z_2 = 0,
\end{equation}
\begin{equation}\label{9.3}
\nabla_{Z_1} \zeta = 0,
\end{equation}
\begin{equation}\label{9.4}
S(Z_1,\zeta) = 0.
\end{equation}

\begin{lem}(\cite{erk})
For a 3-dimensional para-cosymplectic manifold $M^3$, we have
\begin{equation}\label{9.5}
Q Z_1 = \frac{r}{2}[Z_1 - \eta(Z_1)\zeta],
\end{equation}
\begin{equation}\label{9.6}
S(Z_1,Z_2) = \frac{r}{2}[g(Z_1,Z_2)-\eta(Z_1)\eta(Z_2)].
\end{equation}
\end{lem}

\begin{lem}(\cite{erk})
In a para-cosymplectic manifold $M^3$, we get
\begin{equation}\label{9.7}
\zeta r = 0.
\end{equation}
\end{lem}

{\section{\textsf{$h$-ARYS on para-cosymplectic manifolds}}}

Assume that the para-cosymplectic manifold admits an $h$-ARYS $(g,\zeta,\lambda,\alpha,\beta)$. Then (\ref{1.3}) implies
\begin{equation}\label{a.1}
\frac{h}{2}(\pounds_{\zeta} g)(Z_1,Z_2) + \alpha S(Z_1,Z_2) + (\lambda -\frac{\beta}{2}r)g(Z_1,Z_2) = 0,
\end{equation}
which turns into
\begin{equation}\label{a.2}
\frac{h}{2}[g(\nabla_{Z_1} {\zeta},Z_2) + g(Z_1,\nabla_{Z_2} {\zeta})] + \alpha S(Z_1,Z_2) + (\lambda -\frac{\beta}{2}r)g(Z_1,Z_2) = 0.
\end{equation}
Using (\ref{9.3}) in (\ref{a.2}) gives
\begin{equation}
\alpha S(Z_1,Z_2) = -(\lambda-\frac{\beta}{2}r)g(Z_1,Z_2),
\end{equation}
which is an Einstein manifold. Hence we have:\\

\begin{theo}
If a para-cosymplectic manifold admits a proper $h$-ARYS, then the manifold becomes an Einstein manifold.
\end{theo}

\vspace{15cm}

Let $M^3$ admits an $h$-AGRYS. Then from (\ref{2.4}), we get
\begin{equation}\label{b.1}
h\nabla_{Z_1} Df = -\alpha QZ_1 - (\lambda-\frac{\beta}{2}r)Z_1.
\end{equation}
Hence we have
\begin{eqnarray}\label{b.2}
hR(Z_1,Z_2)Df &=& \frac{1}{h}(Z_1 h)[\alpha QZ_2 + (\lambda-\frac{\beta}{2}r)Z_2]\\ \nonumber
&& - \frac{1}{h}(Z_2h)[\alpha QZ_1 + (\lambda-\frac{\beta}{2}r)Z_1]\\ \nonumber
&& -\alpha[(\nabla_{Z_1} Q)Z_2 - (\nabla_{Z_2} Q)Z_1] - (Z_1\lambda)Z_2 + (Z_2\lambda)Z_1\\ \nonumber
&& + \frac{\beta}{2}[(Z_1 r)Z_2 - (Z_2 r)Z_1].
\end{eqnarray}
Using (\ref{9.5}) in (\ref{b.2}) reveals that
\begin{eqnarray}\label{b.3}
hR(Z_1,Z_2)Df &=& \frac{1}{h}(Z_1h)[\alpha QZ_2 + (\lambda-\frac{\beta}{2}r)Z_2]\\ \nonumber
&& - \frac{1}{h}(Z_2 h)[\alpha QZ_1 + (\lambda-\frac{\beta}{2}r)Z_1]\\ \nonumber
&& -\frac{\alpha}{2}(Z_1 r)[Z_2-\eta(Z_2)\zeta] + \frac{\alpha}{2}(Z_2r)[Z_1-\eta(Z_1)\zeta]\\ \nonumber
&& - (Z_1\lambda)Z_2 + (Z_2\lambda)Z_1 + \frac{\beta}{2}[(Z_1 r)Z_2 - (Z_2 r)Z_1].
\end{eqnarray}
If we take $h$ = constant, then the above equation implies
\begin{eqnarray}\label{b.4}
hR(Z_1,Z_2)Df &=&  -\frac{\alpha}{2}(Z_1 r)[Z_2-\eta(Z_2)\zeta] + \frac{\alpha}{2}(Z_2 r)[Z_1-\eta(Z_1)\zeta]\\ \nonumber
&& - (Z_1 \lambda)Z_2 + (Z_2 \lambda)Z_1 + \frac{\beta}{2}[(Z_1 r)Z_2 - (Z_2 r)Z_1].
\end{eqnarray}
Contracting the foregoing equation entails that
\begin{equation}\label{b.5}
hS(Z_2,Df) = (\frac{\alpha}{2}-\beta)Z_2 r + 2(Z_2 \lambda).
\end{equation}
Substituting $Z_1$ by $Df$ in (\ref{9.6}) and equating with (\ref{b.5}), we get
\begin{equation}\label{b.6}
\frac{hr}{2}[Z_2f -(\zeta f)\eta(Z_2)] = (\frac{\alpha}{2}-\beta)Z_2 r + 2(Z_2 \lambda).
\end{equation}
Putting $Z_2 = \zeta$ and using (\ref{9.7}), we infer
\begin{equation}\label{b.7}
\zeta \lambda = 0.
\end{equation}
Taking inner product of (\ref{b.4}) with $\zeta$ and using (\ref{9.1}) gives
\begin{equation}\label{b.8}
-(Z_1\lambda)\eta(Z_2) + (Z_2\lambda)\eta(Z_1) + \frac{\beta}{2}[(Z_1 r)\eta(Z_2) - (Z_2 r)\eta(Z_1)] = 0.
\end{equation}
Setting $Z_2 = \zeta$ in (\ref{b.8}), we get
\begin{equation}\label{b.9}
-(Z_1\lambda) + \frac{\beta}{2}(Z_1 r) = 0.
\end{equation}
If we take $r$ = constant, then (\ref{b.9}) implies
\begin{equation}
Z_1 \lambda = 0,
\end{equation}
which implies $\lambda$ is constant. Therefore we have:\\

\begin{theo}
If a $M^3$ with constant scalar curvature admits an $h$-AGRYS, then the soliton becomes an $h$-gradient Ricci-Yamabe soliton.
\end{theo}

In particular, if we take $\alpha =1$ and $\beta = 0$, then (\ref{b.9}) implies $Z_1\lambda = 0$. Therefore $\lambda$ is constant. Hence we have:
\begin{cor}
An  $h$-almost gradient Ricci soliton  in a $M^3$ becomes an $h$-gradient Ricci soliton.
\end{cor}

\vspace{.9cm}

{\section{\textsf{Examples}}}

{\bf{Example 1.}}
Let us consider $M^3 = \lbrace (x,y,z)\in \mathbb{R}^3 : (x,y,z) \neq (0,0,0)\rbrace$, where $(x,y,z)$ are the standard co-ordinates of $\mathbb{R}^3$.\\
We consider three linearly independent vector fields
\begin{equation} \nonumber
u_1 = e^z\frac{\partial}{\partial x},\hspace{.4cm} u_2 = e^{-z}\frac{\partial}{\partial y},\hspace{.4cm} u_3 = \frac{\partial}{\partial z}.
\end{equation}
Let $g$ be the semi-Riemannian metric defined by
\begin{equation}\nonumber
g(u_1,u_1) = 1,\hspace{.4cm} g(u_2,u_2) = -1,\hspace{.4cm} g(u_3,u_3) = 1,
\end{equation}
\begin{equation}\nonumber
g(u_1,u_2) = 0,\hspace{.4cm} g(u_1,u_3) = 0,\hspace{.4cm} g(u_2,u_3) = 0.
\end{equation}
Let $\eta$ be the 1-form defined by $\eta(Z_1) = g(Z_1,u_3)$ for any $Z_1 \in \chi(M)$.\\
Let $\phi$ be the (1,1)-tensor field defined by
\begin{equation}\nonumber
\phi u_1 = u_2,\hspace{.4cm} \phi u_2 = u_1,\hspace{.4cm} \phi u_3 = 0.
\end{equation}
Using the above relations, we acquire
\begin{equation}\nonumber
\phi^2 Z_1 = Z_1 -\eta(Z_1)u_3,\hspace{.4cm} \eta(u_3) = 1,
\end{equation}
\begin{equation}\nonumber
g(\phi Z_1, \phi Z_2) = -g(Z_1,Z_2) + \eta(Z_1)\eta(Z_2)
\end{equation}
for any $Z_1,Z_2 \in \chi(M)$. Hence for $u_3 = \zeta$, the structure $(\phi,\zeta,\eta,g)$ is an almost paracontact structure on $M$.\\
Using (\ref{6.3}), we have
\begin{equation}\nonumber
\nabla_{u_1} u_1 = -u_3,\hspace{.4cm} \nabla_{u_1} u_2 = 0,\hspace{.4cm} \nabla_{u_1} u_3 = u_1,
\end{equation}
\begin{equation}\nonumber
\nabla_{u_2} u_1 = 0,\hspace{.4cm} \nabla_{u_2} u_2 = u_3,\hspace{.4cm} \nabla_{u_2} u_3 = u_2,
\end{equation}
\begin{equation}\nonumber
\nabla_{u_3} u_1 =0,\hspace{.4cm} \nabla_{u_3} u_2 = 0,\hspace{.4cm} \nabla_{u_3} u_3 =0.
\end{equation}
Hence the manifold is a  para-Kenmotsu manifold. \\
With the help of the above results we can easily obtain
\begin{equation}\nonumber
R(u_1,u_2)u_3 = 0,\hspace{.4cm} R(u_2,u_3)u_3 = - u_2,\hspace{.4cm} R(u_1,u_3)u_3 = -u_1,  
\end{equation}
\begin{equation}\nonumber
R(u_1,u_2)u_2 = u_1,\hspace{.4cm} R(u_2,u_3)u_2 = u_3,\hspace{.4cm} R(u_1,u_3)u_2 = 0,  
\end{equation}
\begin{equation}\nonumber
R(u_1,u_2)u_1 = u_2,\hspace{.4cm} R(u_2,u_3)u_1 = 0,\hspace{.4cm} R(u_1,u_3)u_1 = u_3
\end{equation}
and
\begin{equation}\nonumber
S(u_1,u_1) = -2,\hspace{.4cm} S(u_2,u_2) = 2,\hspace{.4cm} S(u_3,u_3) = -2.
\end{equation}
From the above results we get $r = -6$.\\

Again, suppose that $-f = \lambda = e^z$ and  $2\alpha - 3\beta = 0$. Therefore ${-Df = D\lambda = e^z u_3}$. Hence we get
\begin{eqnarray} \nonumber
&&\nabla_{u_1} \mathcal{D}f = -e^z u_1,\\ \nonumber
&&\nabla_{u_2} \mathcal{D}f = -e^z u_2,\\ \nonumber
&&\nabla_{u_3} \mathcal{D}f =-e^z u_3.
\end{eqnarray}
Therefore, for $2\alpha - 3\beta = 0$ and $h = 1$ ,   equation (\ref{4.1}) is satisfied. Thus $g$ is an $h$-AGRYS with the soliton vector field $V = \mathcal{D}f$, where $-f = \lambda = e^z$ and $2\alpha - 3\beta = 0$. Since $-f = \lambda = e^z$ and $-Df = D\lambda = e^z u_3$, hence Theorem 3.2 is verified. \\\\

{\bf{Example 2.}}
Let $M^3 = \lbrace (x,y,z)\in \mathbb{R}^3 : (x,y,z) \neq (0,0,0)\rbrace$, where $(x,y,z)$ are the standard co-ordinates of $\mathbb{R}^3$.\\
We consider
\begin{equation} \nonumber
v_1 = e^z\frac{\partial}{\partial x} + e^{-z}\frac{\partial}{\partial y},\hspace{.4cm} v_2 = e^z\frac{\partial}{\partial x} - e^{-z}\frac{\partial}{\partial y},\hspace{.4cm} v_3 = \frac{\partial}{\partial z}
\end{equation}
which are linearly independent vector fields. Let the semi-Riemannian metric $g$ be defined by
\begin{equation}\nonumber
g(v_1,v_1) = 1,\hspace{.4cm} g(v_2,v_2) = -1,\hspace{.4cm} g(v_3,v_3) = 1,
\end{equation}
\begin{equation}\nonumber
g(v_1,v_2) = 0,\hspace{.4cm} g(v_1,v_3) = 0,\hspace{.4cm} g(v_2,v_3) = 0.
\end{equation}
The 1-form $\eta$ is defined by $\eta(Z_1) = g(Z_1,v_3)$ for any $Z_1 \in \chi(M)$ and $\phi$ is defined by
\begin{equation}\nonumber
\phi v_1 = v_2,\hspace{.4cm} \phi v_2 = v_1,\hspace{.4cm} \phi v_3 = 0.
\end{equation}
Using the above relations, we acquire
\begin{equation}\nonumber
\phi^2 Z_1 = Z_1 -\eta(Z_1)v_3,\hspace{.4cm} \eta(v_3) = 1,
\end{equation}
\begin{equation}\nonumber
g(\phi Z_1, \phi Z_2) = -g(Z_1,Z_2) + \eta(Z_1)\eta(Z_2)
\end{equation}
for any $Z_1,Z_2 \in \chi(M)$. Hence for $v_3 = \zeta$, the structure $(\phi,\zeta,\eta,g)$ is an almost paracontact structure on $M$.\\
Using (\ref{6.3}), we have
\begin{equation}\nonumber
\nabla_{v_1} v_1 = 0,\hspace{.4cm} \nabla_{v_1} v_2 = v_3,\hspace{.4cm} \nabla_{v_1} v_3 = -v_2,
\end{equation}
\begin{equation}\nonumber
\nabla_{v_2} v_1 = v_3,\hspace{.4cm} \nabla_{v_2} v_2 = 0,\hspace{.4cm} \nabla_{v_2} v_3 = -v_1,
\end{equation}
\begin{equation}\nonumber
\nabla_{v_3} v_1 =0,\hspace{.4cm} \nabla_{v_3} v_2 = 0,\hspace{.4cm} \nabla_{v_3} v_3 =0.
\end{equation}
Hence the manifold is a  para-Sasakian manifold. The components of the curvature tensor and Ricci tensor  are
\begin{equation}\nonumber
R(v_1,v_2)v_3 = 0,\hspace{.4cm} R(v_2,v_3)v_3 = - v_2,\hspace{.4cm} R(v_1,v_3)v_3 = -v_1,
\end{equation}
\begin{equation}\nonumber
R(v_1,v_2)v_2 =v_1,\hspace{.4cm} R(v_2,v_3)v_2 = v_3,\hspace{.4cm} R(v_1,v_3)v_2 = 0,
\end{equation}
\begin{equation}\nonumber
R(v_1,v_2)v_1 = -v_2,\hspace{.4cm} R(v_2,v_3)v_1 = 0,\hspace{.4cm} R(v_1,v_3)v_1 = v_3
\end{equation}
and
\begin{equation}\nonumber
S(v_1,v_1) = -2,\hspace{.4cm} S(v_2,v_2) = 2,\hspace{.4cm} S(v_3,v_3) = -2.
\end{equation}
From the above equations, we obtain $r = -6$.\\

From (\ref{7.3}) we obtain $\alpha S(u_1,u_1)= -3\beta - \lambda$, $\alpha S(u_2,u_2) = 3\beta+\lambda$ and $\alpha S(u_3,u_3) = -3\beta - \lambda$, hence $\lambda = 2\alpha-3\beta$. The data $(g, \zeta, \lambda, \alpha, \beta)$ defines an $h$-ARYS soliton on the para-Sasakian  manifold.\\

\end{document}